# HOW LIKELY IS AN I.I.D. DEGREE SEQUENCE TO BE GRAPHICAL?

By Richard Arratia and Thomas M. Liggett[1]

*University of Southern California and University of California, Los Angeles*

Given i.i.d. positive integer valued random variables $D_1, \ldots, D_n$, one can ask whether there is a simple graph on $n$ vertices so that the degrees of the vertices are $D_1, \ldots, D_n$. We give sufficient conditions on the distribution of $D_i$ for the probability that this be the case to be asymptotically $0$, $\frac{1}{2}$ or strictly between $0$ and $\frac{1}{2}$. These conditions roughly correspond to whether the limit of $nP(D_i \geq n)$ is infinite, zero or strictly positive and finite. This paper is motivated by the problem of modeling large communications networks by random graphs.

**1. Introduction and statement of results.** The growth of the Internet and other large communications systems over the past few decades has led to a significant amount of interest in modeling such systems as random graphs with a large number of vertices. Many books and papers in both mathematics and applied areas such as computer science have resulted from this interest.

Classical mathematical work on random graphs concentrated on the so-called Erdős–Rényi graphs. In these graphs, one takes $n$ vertices and connects pairs of vertices by edges independently with probability $p$ for each pair. While these random graphs have a rich mathematical theory [see Bollobás (2001) or Janson, Luczak and Rucinski (2000), e.g.], it turns out that they are not good choices for modeling many complex networks. To illustrate their limitations, consider the distribution of the degree of a fixed vertex, which in Erdős–Rényi graphs is Binomial with parameters $n-1$ and $p$. If one takes $n$ large and $p$ small in the usual way, this distribution is approximately Poisson. The Poisson distribution has very rapidly decreasing tails.

Received April 2003; revised December 2003.
[1]Supported in part by NSF Grants DMS-00-70465 and DMS-03-01795.

Networks that appear in practice, on the other hand, tend to have degree distributions with very large tails, often lacking a third, or even second, moment. This paper is a contribution to the problem of constructing random graph models in which the distribution of degrees can be prescribed by the modeler.

Many authors have observed that in graphs of complex networks such as the Internet, the degree sequence appears to fit a power law of the form

$$P(D > k) \sim ck^{1-\beta}.$$

One way to mimic the randomness of Erdős–Rényi graphs, while maintaining a power law degree distribution, is to fix a law, sample $n$ times from this law to get a (potential) degree sequence, and then, if there are graphs with this degree sequence, choose a graph uniformly from all graphs with the given degree sequence. Albert and Barabási (2002) provided a good survey of such complex networks and models such as that above. To avoid the uncertainty of sampling, Aiello, Chung and Lu (2001) took a power law and formed a degree sequence by taking $d_i$ (for $1 \leq i \leq n$) to be the (ceiling of the) $(i/n)$th quantile of that law. A survey of rigorous results on graphs with power law degree distributions appeared in Bollobás and Riordan (2003).

To implement the procedure described in the preceding paragraph, two problems must be solved:

1. Determine which degree distributions have the property that if the degrees are chosen by sampling from that distribution, then there is substantial probability that there exists a simple graph with those degrees if $n$ is large.
2. Given that there is such a graph, find efficient algorithms that will make a uniform choice from among these graphs.

This paper is devoted to the first problem. The second problem was treated in several papers: Jerrum and Sinclair (1990), Jerrum, McKay and Sinclair (1992) and Steger and Wormald (1999). Note that while the titles of these papers often suggest that they treat the problem only for regular graphs (in which the degrees are all the same), the algorithms proposed usually apply to nonregular graphs as well, provided that the degrees do not vary too much. In fact, conditions on the degree sequence are given which imply that the algorithm works well, and then it is observed that these conditions are satisfied by regular graphs.

We now turn to the mathematical formulation of our problem and to the statement of our main result. Consider a sequence of positive integers $d_1, \ldots, d_n$. The decreasing rearrangement of this sequence is $m_1, \ldots, m_n$, where $\{m_1, \ldots, m_n\} = \{d_1, \ldots, d_n\}$ as multisets and $m_1 \geq m_2 \geq \cdots \geq m_n$. Erdős and Gallai (1960) showed that there exists a simple graph (i.e., one with no loops or multiple edges) with $n$ vertices whose respective degrees are $d_1, \ldots, d_n$ if and only if $\sum_i d_i$ is even and their decreasing rearrangement satisfies

$$(1.1) \qquad \sum_{i=1}^{j} m_i \leq j(j-1) + \sum_{i=j+1}^{n} \min(j, m_i), \qquad 1 \leq j \leq n.$$

Sequences with these properties are called graphical. Some equivalent conditions for a sequence to be graphical can be found in Sierksma and Hoogeveen (1991).

The necessity of (1.1) is quite easy to check. The left-hand side of (1.1) is the total degree of the first $j$ vertices. The first term on the right-hand side is the maximum total degree of these vertices coming from edges that connect vertices in $\{1, \ldots, j\}$ to other vertices in this set. The second term is a bound on the total degree of the vertices in $\{1, \ldots, j\}$ coming from edges that connect to vertices in $\{j+1, \ldots, n\}$. The proof of the sufficiency of (1.1) is more difficult. Proofs can be found in several papers, including Sierksma and Hoogeveen (1991) and Choudum (1986). The latter paper provides a simple inductive construction of the required graph. Briefly, suppose

$$m_1 = \cdots = m_j > m_{j+1} \geq \cdots \geq m_n \geq 1$$

satisfies (1.1). While not obvious, it can be checked that the sequence of degrees

$$m_1, \ldots, m_{j-1}, m_j - 1, m_{j+1}, \ldots, m_{n-1}, m_n - 1$$

also satisfies (1.1) (and has smaller total degree). A graph with the original degree sequence can be constructed from one with the modified degree sequence either by adding one edge or by deleting two edges and adding three edges, according to whether the edge $(j, n)$ is already in the graph with smaller total degree.

To state our main result, let $D$ be a positive integer valued random variable and let $D_1, \ldots, D_n$ be an i.i.d. sequence of random variables with the distribution of $D$. Let $\Lambda$ be the set of limit points of the sequence of probabilities

$$P((D_1, \ldots, D_n) \text{ is graphical}).$$

If $0 < P(D \text{ is even }) < 1$, then the parity restriction forces $\Lambda \subset [0, \frac{1}{2}]$, as we see below. We are interested in finding sufficient conditions for $\Lambda = \{0\}, \Lambda = \{\frac{1}{2}\}$, or $\Lambda \subset (0, \frac{1}{2})$. Our main result is given below. To simplify its statement, we consider only the case in which $D$ takes both even and odd values with positive probability.

THEOREM. *Suppose that $0 < P(D \text{ is even}) < 1$. Then the following statements hold*:

(a) If $\lim_{n\to\infty} nP(D \geq n) = \infty$, then $\Lambda = \{0\}$.
(b) If $\lim_{n\to\infty} nP(D \geq n) = c$, where $0 < c < \infty$, and

$$\sum_n \left[\frac{c}{n} - P(D \geq n)\right] = \infty, \tag{1.2}$$

then $\Lambda = \{0\}$.
(c) If $\lim_{n\to\infty} nP(D \geq n) = c$, where $0 < c < \infty$,

$$0 < \liminf_{n\to\infty} n^2 P(D = n) \leq \limsup_{n\to\infty} n^2 P(D = n) < \infty \tag{1.3}$$

and

$$\sum_n \left|\frac{c}{n} - P(D \geq n)\right| < \infty, \tag{1.4}$$

then $\Lambda \subset (0, \frac{1}{2})$.
(d) If $ED < \infty$ or if

$$\sup_n n(\log n) P(D \geq n) < \infty, \tag{1.5}$$

then $\Lambda = \{\frac{1}{2}\}$.

Note that it follows from parts (b) and (c) of the Theorem that $\Lambda$ is not monotone in the distribution of $D$. The main interest in part (d) is the case in which $ED < \infty$. We include the result when (1.5) holds primarily because this shows that $ED < \infty$ is not necessary and sufficient for $\Lambda = \{\frac{1}{2}\}$.

Another version of this problem was considered earlier. To describe it, recall that a partition of a positive integer $N$ is a sequence of integers $\gamma_1 \geq \gamma_2 \geq \cdots \geq \gamma_n > 0$ so that

$$N = \sum_{i=1}^n \gamma_i.$$

In 1982, Wilf conjectured that if $\Gamma_N$ is chosen at random uniformly from all partitions of $N$, then

$$\lim_{N\to\infty} P(\Gamma_N \text{ is graphical}) = 0.$$

This conjecture was proved by Pittel (1999).

**2. Proof of the Theorem.** This section is devoted to the proof of the Theorem. The order in which the four parts are stated and proved is according to the size of the tails of the distribution of $D$—from long tails to short. However, the proof of part (c) is much longer than the proofs of the other parts, so you might want to read the proofs in the order (a), (b), (d) and (c).

Before entering into the main part of the proof, we deal with the easy parity issue. If $P(D \text{ is even}) = 1$, then clearly $P(D_1 + \cdots + D_n \text{ is even}) = 1$ as well. Similarly, if $P(D \text{ is odd}) = 1$, then $P(D_1 + \cdots + D_n \text{ is even}) = 0$ or 1 according to whether $n$ is odd or even.

PROPOSITION. *If $0 < P(D \text{ is even}) < 1$, then*

$$\lim_{n \to \infty} P(D_1 + \cdots + D_n \text{ is even}) = \tfrac{1}{2}.$$

PROOF. The result follows immediately from

$$2P(D_1 + \cdots + D_n \text{ is even}) - 1 = E((-1)^{D_1 + \cdots + D_n}) = [E(-1)^D]^n$$

and the fact that $|E(-1)^D| < 1$. □

We turn now to the proof of the Theorem, for which it is useful to restate (1.1) in a more convenient form. Let $M_n(1) \geq M_n(2) \geq \cdots \geq M_n(n)$ be the decreasing rearrangement of $D_1, \ldots, D_n$. Then the inequality in (1.1) becomes

$$(2.1) \qquad \sum_{i=1}^{j} M_n(i) \leq j(j-1) + \sum_{i=j+1}^{n} \min(j, M_n(i)),$$

which can be written as

$$(2.2) \qquad \sum_{i=1}^{j} M_n(i) \leq j(j-1) + \sum_{1 \leq l \leq j < i \leq n} \mathbb{1}_{\{M_n(i) \geq l\}},$$

where $\mathbb{1}_A$ denotes the indicator of the set $A$.

PROOF OF PART (a). Considering only the case $j = 1$ in (2.2), we see that

$$P((D_1, \ldots, D_n) \text{ is graphical}) \leq P(M_n(1) \leq n - 1) = [1 - P(D \geq n)]^n,$$

which tends to 0 as $n \to \infty$ by assumption. □

The proofs of parts (b)–(d) of the Theorem are based on a useful representation of $M_n(1), \ldots, M_n(n)$. To describe it, let $X_1, X_2, \ldots$ be independent unit exponentials, and let $g : [0, \infty) \to \{1, 2, \ldots\}$ be a right continuous, increasing function. Take $D_i = g(1/X_i)$, so that

$$(2.3) \qquad P(D \geq n) = 1 - \exp\left(-\frac{1}{g^{-1}(n)}\right),$$

where $g^{-1}(n)$ is defined by $g(u) \geq n$ iff $u \geq g^{-1}(n)$. Note that the distribution of $D_i$ uniquely determines the values $g^{-1}(n)$, and hence uniquely determines the function $g$, since $g$ is taken to be increasing and right continuous. This construction of the $D_i$'s is a variant of the familiar quantile construction of an arbitrary random variable $X$ as $X = F^{-1}(U)$ in terms of a random variable $U$ that is uniform on $[0,1]$.

For fixed $n$, let $X_{(1)}, \ldots, X_{(n)}$ be the increasing rearrangement of $X_1, \ldots, X_n$. Then

$$M_n(i) = g\left(\frac{1}{X_{(i)}}\right), \qquad 1 \leq i \leq n.$$

Elementary properties of the exponential distribution imply that

$$(X_{(1)}, \ldots, X_{(n)}) \stackrel{d}{=} \left(\frac{X_1}{n}, \frac{X_1}{n} + \frac{X_2}{n-1}, \ldots, \frac{X_1}{n} + \frac{X_2}{n-1} + \cdots + \frac{X_n}{1}\right),$$

where $\stackrel{d}{=}$ denotes equality in distribution. [This useful representation was first observed and exploited by Epstein and Sobel (1953) and Rényi (1953).] It follows from this that

$$(M_n(1), M_n(2), \ldots, M_n(n))$$

(2.4)
$$\stackrel{d}{=} \left(g\left(\left(\frac{X_1}{n}\right)^{-1}\right), g\left(\left(\frac{X_1}{n} + \frac{X_2}{n-1}\right)^{-1}\right), \ldots, \right.$$
$$\left. g\left(\left(\frac{X_1}{n} + \frac{X_2}{n-1} + \cdots + \frac{X_n}{1}\right)^{-1}\right)\right).$$

As usual, we let $S_n = X_1 + \cdots + X_n$ denote the partial sums of the $X_i$'s.

PROOF OF PART (b). Since $nP(D \geq n) \to c$, (2.3) implies that $g^{-1}(n) \sim n/c$ and hence that

(2.5) $$\lim_{x \to \infty} \frac{g(x)}{x} = c.$$

By (2.4) and (2.5),

(2.6) $$\left(\frac{M_n(1)}{n}, \ldots, \frac{M_n(j)}{n}\right) \stackrel{d}{\to} \left(\frac{c}{S_1}, \ldots, \frac{c}{S_j}\right)$$

for $j \geq 1$, where $\stackrel{d}{\to}$ denotes convergence in distribution. Statements of this type are well known, of course; see, for example, Section 2.3 of Leadbetter, Lindgren and Rootzén (1983). However, especially in the proof of part (c) of the Theorem, we need much more detailed information about the asymptotic behavior of $(M_n(1), \ldots, M_n(n))$ than such results provide.

By the strong law of large numbers,

$$\lim_{n \to \infty} \frac{1}{n} \sum_{i=j+1}^{n} \mathbb{1}_{\{M_n(i) \geq l\}} = \lim_{n \to \infty} \frac{1}{n} \sum_{i=1}^{n} \mathbb{1}_{\{M_n(i) \geq l\}}$$

$$= \lim_{n \to \infty} \frac{1}{n} \sum_{i=1}^{n} \mathbb{1}_{\{D_i \geq l\}} = P(D \geq l) \quad \text{a.s.}$$

for each $j$. Dividing (2.2) by $n$ and passing to the limit as $n \to \infty$, it follows that

$$\limsup_{n \to \infty} P((D_1, \ldots, D_n) \text{ is graphical}) \leq P\left(\sum_{i=1}^{j} \frac{c}{S_i} \leq \sum_{i=1}^{j} P(D \geq i) \text{ for all } j \geq 1\right).$$

So, it suffices to show that

$$\sum_{i=1}^{\infty} \left(\frac{c}{S_i} - P(D \geq i)\right) = \infty \quad \text{a.s.}$$

This follows from (1.2), since

$$\sum_{i=1}^{\infty} \left(\frac{c}{S_i} - \frac{c}{i}\right) = c \sum_{i=1}^{\infty} \frac{i - S_i}{i S_i},$$

which converges absolutely a.s. by the law of the iterated logarithm applied to the sequence $S_1, S_2, \ldots$. □

PROOF OF PART (c). We begin by observing that (1.3) implies that there is a constant $C$ so that

(2.7) $$|g(u) - g(v)| \leq C|u - v| + 1$$

for all $u, v \geq 0$. To see this, use (1.3) and (2.3) to conclude that there is an $\varepsilon > 0$ and an $N$ so that

$$g^{-1}(n+1) - g^{-1}(n) \geq \varepsilon$$

for all $n \geq N$. Now take $0 \leq u \leq v$ and let $m = g(u), n = g(v)$. Without loss of generality, we may assume $N \leq m < n$ since (2.7) is automatic if $n = m$ and once we know (2.7) holds for large $u, v$, it can be made to hold for other $u, v$ by increasing the value of $C$. In this case,

$$g^{-1}(m) \leq u < g^{-1}(m+1) \quad \text{and} \quad g^{-1}(n) \leq v < g^{-1}(n+1),$$

so that

$$v - u \geq g^{-1}(n) - g^{-1}(m+1) \geq \varepsilon(n - m - 1).$$

This gives (2.7) with $C = \varepsilon^{-1}$.

In what follows, we use $O(1)$ in distributional statements to represent possibly random quantities that are bounded by a universal constant. In limiting statements of the form $Y_n \xrightarrow{d} Y + O(1)$, the meaning is that the sequence of distributions of $Y_1, Y_2, \ldots$ is tight and that all weak limits of the distribution of $Y_n$ lie between the distribution of $Y \pm$ a universal constant.

Since the proof of part (c) is quite long, we break it up into several lemmas. The first determines the asymptotics of the distribution of the left-hand side of (2.1), when centered appropriately. The second lemma carries out the additional analysis that is needed to deal with the right-hand side of (2.1). The final lemma provides bounds for the nonrandom quantities that arise in the centering in the first two lemmas. $\square$

LEMMA 1. *For every $N \geq 1$,*

$$\max_{N \leq j \leq n} \sum_{i=1}^{j} \left( \frac{M_n(i)}{n} - \frac{1}{n} g\left( \left( \sum_{l=1}^{i} \frac{1}{n-l+1} \right)^{-1} \right) \right)$$

$$\xrightarrow{d} c \sup_{j \geq N} \left( \sum_{i=1}^{j} \left( \frac{1}{S_i} - \frac{1}{i} \right) \right) + O(1).$$

PROOF. In the proof, we need to treat the summands above separately in three regimes: $i \leq \frac{n}{2}, \frac{n}{2} \leq i \leq n - n^\theta$ and $n - n^\theta \leq i \leq n$. We begin with some estimates that are needed in all three cases. Define

$$Z_n(i) = \frac{M_n(i)}{n} - \frac{1}{n} g\left( \left( \sum_{l=1}^{i} \frac{1}{n-l+1} \right)^{-1} \right)$$

for $1 \leq i \leq n$. Then jointly in $i$,

(2.8)
$$n|Z_n(i)| \stackrel{d}{=} \left| g\left( \left( \sum_{l=1}^{i} \frac{X_l}{n-l+1} \right)^{-1} \right) - g\left( \left( \sum_{l=1}^{i} \frac{1}{n-l+1} \right)^{-1} \right) \right|$$

$$\stackrel{d}{\leq} \left( C \left| \sum_{l=1}^{i} \frac{X_l - 1}{n-l+1} \right| \right) \Big/ \left( \sum_{l=1}^{i} \frac{X_l}{n-l+1} \sum_{l=1}^{i} \frac{1}{n-l+1} \right) + O(1)$$

by (2.4) and (2.7). To bound the right-hand side of (2.8), sum by parts to get

$$\sum_{l=1}^{i} \frac{X_l - 1}{n - l + 1} = \sum_{l=1}^{i} \frac{(S_l - l) - (S_{l-1} - l + 1)}{n - l + 1}$$

$$= \frac{S_i - i}{n - i + 1} - \sum_{l=1}^{i-1} \frac{S_l - l}{(n-l)(n-l+1)},$$

where again $S_k = X_1 + \cdots + X_k$. Therefore, letting

$$S = \sup_{k \geq 1} \frac{|S_k - k|}{\sqrt{k \log(k+1)}},$$

which is finite a.s. by the law of the iterated logarithm, we see that for $i \leq n$,

$$\text{(2.9)} \qquad \left| \sum_{l=1}^{i} \frac{X_l - 1}{n - l + 1} \right| \leq 2S \frac{\sqrt{i \log(i+1)}}{n - i + 1}.$$

Combining (2.8) and (2.9) gives

$$\text{(2.10)} \quad n|Z_n(i)| \overset{d}{\leq} 2CS\sqrt{i \log(i+1)}$$
$$\times \left( \left( \sum_{l=1}^{i} \frac{n-i+1}{n-l+1} - 2S\sqrt{i \log(i+1)} \right) \sum_{l=1}^{i} \frac{1}{n-l+1} \right)^{-1} + O(1),$$

provided that the denominator is positive. Here again, the distributional inequality holds jointly for $1 \leq i \leq n$.

For $1 \leq i \leq \frac{n}{2}$, we may replace the $n - l + 1$ by $n$ and the $n - i + 1$ by $\frac{n}{2}$ on the right-hand side of (2.10) and divide by $n$ to get

$$\text{(2.11)} \qquad |Z_n(i)| \overset{d}{\leq} \frac{4CS\sqrt{i \log(i+1)}}{(i - 4S\sqrt{i \log(i+1)})i} + \frac{1}{n}O(1),$$

again provided that the denominator on the right-hand side is positive.

By (2.5) and (2.6),

$$Z_n(i) \overset{d}{\to} c\left( \frac{1}{S_i} - \frac{1}{i} \right)$$

jointly for any finite number of $i$'s. Therefore,

$$\max_{N \leq j \leq N+m} \sum_{i=1}^{j} Z_n(i) \overset{d}{\to} c \max_{N \leq j \leq N+m} \left( \sum_{i=1}^{j} \left( \frac{1}{S_i} - \frac{1}{i} \right) \right)$$

for every $N, m \geq 1$. We wish to have a version of this statement that corresponds to $m = \infty$; for this, we need some domination. This is provided by (2.11) as we now argue.

Note first that the denominator of the first term on the right-hand side of (2.11) is positive for sufficiently large $i$ (depending on the value of $S$). Therefore, since

$$\sum_{i=1}^{\infty} \frac{\sqrt{i \log(i+1)}}{i^2} < \infty,$$

it follows that there is a function $f(s)$ so that $i \geq f(s)$ implies $i - 4s\sqrt{i \log(i+1)} > 0$ and

$$\sum_{i \geq f(s)} \frac{4s\sqrt{i \log(i+1)}}{(i - 4s\sqrt{i \log(i+1)})i} \leq 1, \qquad s \geq 0.$$

Combining this with (2.11), it follows that there is a constant $K$ [that can be taken to be the sum of the constant in the $O(1)$ term in (2.11) and $C$] so that

$$P(N+m \geq f(S)) \leq P\left(\sum_{N+m<i\leq n/2} |Z_n(i)| \leq K\right).$$

The final ingredient is the inequality

$$\max_{N\leq j\leq N+m} \sum_{i=1}^{j} Z_n(i) \leq \max_{N\leq j\leq n/2} \sum_{i=1}^{j} Z_n(i)$$

$$\leq \max_{N\leq j\leq N+m} \sum_{i=1}^{j} Z_n(i) + \sum_{N+m<i\leq n/2} |Z_n(i)|$$

for $N+m \leq \frac{n}{2}$. Combining these observations, the conclusion is that for any $N \geq 1$,

$$(2.12) \qquad \max_{N\leq j\leq n/2} \sum_{i=1}^{j} Z_n(i) \xrightarrow{d} c \sup_{j\geq N} \left(\sum_{i=1}^{j}\left(\frac{1}{S_i} - \frac{1}{i}\right)\right) + O(1).$$

Next, we consider terms that correspond to $n - n^\theta \leq i \leq n$, where $0 < \theta < 1$. For such an $i$, the right-hand side of (2.8) is at most

$$(2.13) \qquad \left(C\left(\sum_{l=1}^{n} \frac{X_l + 1}{n-l+1}\right)\right) \bigg/ \left(\sum_{l=1}^{\lfloor n-n^\theta \rfloor} \frac{X_l}{n-l+1} \sum_{l=1}^{\lfloor n-n^\theta \rfloor} \frac{1}{n-l+1}\right) + O(1),$$

where $\lfloor \cdot \rfloor$ is the greatest integer function. Since

$$\sum_{l=1}^{n} \frac{X_l}{n-l+1}$$

has mean that is asymptotic to $\log n$ and variance that is bounded in $n$,

$$\frac{1}{\log n} \sum_{l=1}^{n} \frac{X_l}{n-l+1} \to 1$$

in probability. Similiarly,

$$\frac{1}{\log n} \sum_{l=1}^{\lfloor n-n^\theta \rfloor} \frac{X_l}{n-l+1} \to 1-\theta$$

in probability. Therefore (2.13) is asymptotic to

$$\frac{2C}{(1-\theta)^2 \log n}.$$

It follows that

$$\text{(2.14)} \quad \sum_{n-n^\theta \leq i \leq n} \left| \frac{M_n(i)}{n} - \frac{1}{n} g\left(\left(\sum_{l=1}^{i} \frac{1}{n-l+1}\right)^{-1}\right) \right| \xrightarrow{d} O(1).$$

Finally, we consider the intermediate terms $\frac{n}{2} \leq i \leq n - n^\theta$. In this case,

$$\sum_{l=1}^{i} \frac{1}{n-l+1} \geq \log \frac{n+1}{n-i+1} \geq \log \frac{2n+2}{n+2}$$

and $n - i + 1 \geq n^\theta$, so that the right-hand side of (2.10) is at most

$$\text{(2.15)} \quad (2CS\sqrt{n \log(n+1)}) \Big/ \left( \left( n^\theta \log \frac{2n+2}{n+2} \right. \right.$$
$$\left. \left. - 2S\sqrt{n \log(n+1)} \right) \log \frac{2n+2}{n+2} \right) + O(1),$$

provided that

$$\text{(2.16)} \quad n^\theta \log \frac{2n+2}{n+2} > 2S\sqrt{n \log(n+1)}.$$

If $\frac{1}{2} < \theta < 1$, (2.16) holds eventually a.s. and then (2.15) tends to zero as $n \to \infty$. It follows that

$$\text{(2.17)} \quad \sum_{n/2 \leq i \leq n-n^\theta} \left| \frac{M_n(i)}{n} - \frac{1}{n} g\left(\left(\sum_{l=1}^{i} \frac{1}{n-l+1}\right)^{-1}\right) \right| \xrightarrow{d} O(1)$$

in probability. Combining (2.12), (2.14) and (2.17), we see that for every $N \geq 1$,

$$\text{(2.18)} \quad \max_{N \leq j \leq n} \sum_{i=1}^{j} \left( \frac{M_n(i)}{n} - \frac{1}{n} g\left(\left(\sum_{l=1}^{i} \frac{1}{n-l+1}\right)^{-1}\right) \right)$$
$$\xrightarrow{d} c \sup_{j \geq N} \left( \sum_{i=1}^{j} \left( \frac{1}{S_i} - \frac{1}{i} \right) \right) + O(1),$$

as required. □

To prepare for the second lemma, define

$$\gamma(j, n) = \sum_{i=1}^{j} g\left(\left(\sum_{l=1}^{i} \frac{1}{n-l+1}\right)^{-1}\right)$$
$$+ \sum_{i=1}^{j} g\left(\min\left(j, \left(\sum_{l=1}^{i} \frac{1}{n-l+1}\right)^{-1}\right)\right) - nE \min(j, D).$$

LEMMA 2. *For every $N \geq 1$,*

$$\max_{N \leq j \leq n} \frac{1}{n}\left[\sum_{i=1}^{j} M_n(i) - j(j-1) - \sum_{i=j+1}^{n} \min(j, M_n(i)) - \gamma(j,n)\right]$$

$$\xrightarrow{d} c \sup_{j \geq N}\left(\sum_{i=1}^{j}\left(\frac{1}{S_i} - \frac{1}{i}\right)\right) + O(1).$$

PROOF. Write the right-hand side of (2.1) as

(2.19) $$j(j-1) + \sum_{i=1}^{n} \min(j, D_i) - \sum_{i=1}^{j} \min(j, M_n(i)).$$

The last term is easy to deal with. Since $|\min(j,u) - \min(j,v)| \leq |u-v|$, the estimates in (2.8), (2.11), (2.14) and (2.17) apply just as well when $M_n(i)$ is replaced by $\min(j, M_n(i))$ and

$$g\left(\left(\sum_{l=1}^{i} \frac{1}{n-l+1}\right)^{-1}\right) \quad \text{is replaced by} \quad g\left(\min\left(j, \left(\sum_{l=1}^{i} \frac{1}{n-l+1}\right)^{-1}\right)\right).$$

Now, however, (2.6) is replaced by

$$\left(\frac{\min(j, M_n(1))}{n}, \frac{\min(j, M_n(2))}{n}, \ldots, \frac{\min(j, M_n(j))}{n}\right) \xrightarrow{d} (0, 0, \ldots, 0).$$

Therefore, (2.18) is replaced by

(2.20) $$\max_{1 \leq j \leq n} \sum_{i=1}^{j}\left|\frac{\min(j, M_n(i))}{n} - \frac{1}{n}g\left(\min\left(j, \left(\sum_{l=1}^{i} \frac{1}{n-l+1}\right)^{-1}\right)\right)\right| \xrightarrow{d} O(1).$$

To handle the first two terms in (2.19), we need the asymptotics for moments of $\min(j, D)$ that follow from $nP(D \geq n) \to c$:

$$E\min(j, D) = \sum_{l=1}^{j} P(D \geq l) \sim c \log j,$$

$$E[\min(j, D)]^2 = \sum_{l=1}^{j}(2l-1)P(D \geq l) \sim 2cj,$$

$$E[\min(j, D)]^3 = \sum_{l=1}^{j}(3l^2 - 3l + 1)P(D \geq l) \sim \tfrac{3}{2}cj^2,$$

$$E[\min(j, D)]^4 = \sum_{l=1}^{j}(4l^3 - 6l^2 + 4l - 1)P(D \geq l) \sim \tfrac{4}{3}cj^3.$$

It follows that there is a $C$ so that for $1 \leq j \leq n$,

$$E\left[\frac{1}{n}\sum_{i=1}^{n}[\min(j, D_i) - E\min(j, D)]\right]^4 \leq C\frac{j^2}{n^2}.$$

Now take $\frac{1}{2} < \theta < \frac{2}{3}$. Then

$$E \max_{1 \leq j \leq n^\theta}\left[\frac{1}{n}\sum_{i=1}^{n}[\min(j, D_i) - E\min(j, D)]\right]^4 \leq C \sum_{1 \leq j \leq n^\theta} \frac{j^2}{n^2} \to 0,$$

so that

(2.21) $$\max_{1 \leq j \leq n^\theta}\left|\frac{1}{n}\sum_{i=1}^{n}[\min(j, D_i) - E\min(j, D)]\right| \xrightarrow{d} 0.$$

On the other hand,

$$\min_{n^\theta \leq j \leq n}\frac{1}{n}\left[j(j-1) + \sum_{i=1}^{n}[\min(j, D_i) - E\min(j, D)]\right]$$
$$\geq \frac{n^\theta(n^\theta - 1)}{n} - E\min(n, D),$$

which tends to $\infty$. Combining this with (2.6), (2.11), (2.14), (2.17), (2.20) and (2.21), we conclude that for $N \geq 1$,

(2.22) $$\max_{N \leq j \leq n}\frac{1}{n}\left[\sum_{i=1}^{j}M_n(i) - j(j-1) - \sum_{i=j+1}^{n}\min(j, M_n(i)) - \gamma(j, n)\right]$$
$$\xrightarrow{d} c \sup_{j \geq N}\left(\sum_{i=1}^{j}\left(\frac{1}{S_i} - \frac{1}{i}\right)\right) + O(1),$$

as required. $\square$

Finally, we need to bound the $\gamma(j, n)$'s below.

LEMMA 3. *There exists*

(2.23) $$\inf_{n \geq 1}\frac{1}{n}\min_{1 \leq j \leq n}\gamma(j, n) > -\infty.$$

PROOF. Since we need only a lower bound and $g$ is nonnegative, we may neglect the middle term in the definition of $\gamma(j, n)$. So, using (1.4), it sufficient to show that

$$\sum_{i=1}^{j}\left[\frac{1}{n}g\left(\left(\sum_{l=1}^{i}\frac{1}{n-l+1}\right)^{-1}\right) - \frac{c}{i}\right]$$

is bounded below.

To do so, note that
$$\sum_{l=1}^{i} \frac{1}{n-l+1} \leq \log \frac{n}{n-i},$$
so that
$$\frac{1}{n} g\left(\left(\sum_{l=1}^{i} \frac{1}{n-l+1}\right)^{-1}\right) - c\log\left(1 + \frac{1}{i}\right)$$
$$\geq \int_{i/n}^{(i+1)/n} \left[g\left(\frac{1}{|\log(1-x)|}\right) - \frac{c}{x}\right] dx.$$

Hence, it suffices to prove that
$$\int_0^1 \left|g\left(\frac{1}{|\log(1-x)|}\right) - \frac{c}{x}\right| dx < \infty.$$

The only place where this integrability is an issue is at $x = 0$. Making a change of variables, we see that it suffices to check
$$\int_1^\infty \frac{|g(y) - cy|}{y^2} dy < \infty.$$

So, it is enough to check that
$$\sum_n \int_{g^{-1}(n)}^{g^{-1}(n+1)} \frac{|n - cy|}{y^2} dy < \infty$$
or that
$$\sum_n \frac{\max(|n - cg^{-1}(n)|, |n - cg^{-1}(n+1)|)(g^{-1}(n+1) - g^{-1}(n))}{(g^{-1}(n))^2} < \infty.$$

This follows from (1.3), (1.4) and (2.3). □

We are now in a position to complete the proof of part (c) of the theorem using the three lemmas. Combining (2.22) and (2.23), we see that there is a constant $K$ so that
$$\liminf_{n \to \infty} P\left(\max_{N \leq j \leq n} \frac{1}{n}\left[\sum_{i=1}^{j} M_n(i) - j(j-1) - \sum_{i=j+1}^{n} \min(j, M_n(i))\right] \leq 0\right)$$
$$\geq P\left(c \sup_{j \geq N}\left(\sum_{i=1}^{j}\left(\frac{1}{S_i} - \frac{1}{i}\right)\right) \leq -K\right).$$

Note that the distribution of
$$\sup_{j \geq N} \sum_{i=1}^{j} \left( \frac{1}{S_i} - \frac{1}{i} \right)$$
has support
$$\left[ -\sum_{i=1}^{N} \frac{1}{i}, \infty \right).$$
Therefore $N$ can be chosen large enough that the probability on the right-hand side above is strictly positive. By (2.6), for each $j$,
$$\frac{1}{n} \left[ \sum_{i=1}^{j} M_n(i) - \sum_{i=j+1}^{n} \min(j, M_n(i)) \right] \xrightarrow{d} \sum_{i=1}^{j} \left[ \frac{c}{S_i} - P(D \geq i) \right],$$
which has positive probability of being negative. Since all these statements hold jointly in $j$, it follows that $\liminf_n P(A_n) > 0$, where
$$A_n = \left\{ \sum_{i=1}^{j} M_n(i) \leq j(j-1) + \sum_{i=j+1}^{n} \min(j, M_n(i)) \text{ for all } 1 \leq j \leq n \right\}.$$
On the other hand,
$$P(A_n) \leq P(M_n(1) \leq n - 1) \to e^{-c} < 1.$$
It follows that all limit points of $P(A_n)$ lie in $(0, 1)$.

To complete the proof of part (c), we need to show that $A_n$ and the event
$$\{D_1 + \cdots + D_n \text{ is even}\}$$
are asymptotically independent. To do so, fix $m$ and let $\mathcal{G}$ be the $\sigma$-algebra generated by $\{D_k, k > m\}$, let $A$ be any event, let $B$ be any event depending on $\{D_1, \ldots, D_m\}$ and let $C$ be any event in $\mathcal{G}$. Then
$$P(A \cap B \cap C) - \tfrac{1}{2} P(A \cap C)$$
$$= E[\mathbb{1}_A - E(A|\mathcal{G}), B \cap C] - \tfrac{1}{2} E[\mathbb{1}_A - E(A|\mathcal{G}), C]$$
$$+ P(A \cap C)[P(B) - \tfrac{1}{2}].$$
Therefore,
$$|P(A \cap B \cap C) - \tfrac{1}{2} P(A \cap C)| \leq 2 E|\mathbb{1}_A - E(A|\mathcal{G})| + |P(B) - \tfrac{1}{2}|.$$
Applying this inequality to $B = \{D_1 + \cdots + D_m \text{ is even}\}$, $C = \{D_{m+1} + \cdots + D_n \text{ is even}\}$ and their complements, we have the inequality
$$|P(A, D_1 + \cdots + D_n \text{ is even}) - \tfrac{1}{2} P(A)|$$
$$= |P(A \cap B \cap C) + P(A \cap B^c \cap C^c) - \tfrac{1}{2} P(A \cap C) - \tfrac{1}{2} P(A \cap C^c)|$$
$$\leq 4 E|\mathbb{1}_A - E(A|\mathcal{G})| + |2P(D_1 + \cdots + D_m \text{ is even}) - 1|.$$

Applying this to the $A_n$'s above and using the Proposition, we see that to prove $\Lambda \subset (0, \frac{1}{2})$, it suffices to prove
$$\lim_{n \to \infty} E|\mathbb{1}_{A_n} - P(A_n|\mathcal{G})| = 0.$$
To do this, let $\widetilde{D}_1, \ldots, \widetilde{D}_m$ be independent, independent of $\{D_1, D_2, \ldots\}$ and have the same distribution as $D_i$. Quantities defined in terms of the sequence
$$\{\widetilde{D}_1, \ldots, \widetilde{D}_m, D_{m+1}, \ldots\}$$
are denoted by the earlier notation with a tilde. Then it suffices to prove that
$$\lim_{n \to \infty} P(A_n \Delta \widetilde{A}_n) = 0,$$
where $\Delta$ denotes the symmetric difference. To check that this is enough, note that since $A_n$ and $\widetilde{A}_n$ are conditionally independent given $\mathcal{G}$, we have
$$[E|\mathbb{1}_{A_n} - P(A_n|\mathcal{G})|]^2 \leq E[\mathbb{1}_{A_n} - P(A_n|\mathcal{G})]^2$$
$$= \tfrac{1}{2} E[\mathbb{1}_{A_n} - P(A_n|\mathcal{G}) - \mathbb{1}_{\widetilde{A}_n} + P(\widetilde{A}_n|\mathcal{G})]^2$$
$$\leq E[\mathbb{1}_{A_n} - \mathbb{1}_{\widetilde{A}_n}]^2 \leq 2P(A_n \Delta \widetilde{A}_n).$$

Unfortunately, it is not too easy to see that the right-hand side above tends to zero as required. However, we can apply these same considerations to certain events that are larger (resp. smaller) than these $A_n$'s to complete the proof. For the larger ones, simply use $A_n = \{M_n(1) \leq n - 1\}$. With this choice, it is clear that $P(A_n \Delta \widetilde{A}_n) \to 0$. For the smaller ones, use $A_n = \{W_n \leq -\varepsilon\}$ for $\varepsilon$ small and positive, where
$$W_n = \max_{1 \leq j \leq n} \frac{1}{n}\left[\sum_{i=1}^{j} M_n(i) - j(j-1) - \sum_{i=j+1}^{n} \min(j, M_n(i))\right].$$
Then $W_n - \widetilde{W}_n \to 0$ in probability, since
$$\sum_{i=1}^{n} |M_n(i) - \widetilde{M}_n(i)| \leq \sum_{i=1}^{m} |D_i - \widetilde{D}_i|.$$
Therefore, passing to subsequences so that the distributional limit of $W_n$ [which might have some mass at $-\infty$ but has none at $+\infty$ by (2.22) and (2.23)], we see that $P(A_n \Delta \widetilde{A}_n) \to 0$ provided that this distributional limit has no mass at $\varepsilon$. Now choose such an $\varepsilon$.

PROOF OF PART (d).  By (2.1) and the Proposition, it suffices to show that

(2.24) $$\max_{1 \leq j \leq n}\left[2\sum_{i=1}^{j} M_n(i) - j(j-1) - \sum_{i=1}^{n} \min(j, D_i)\right] \to -\infty$$

in probability. We consider various ranges of $j$'s separately. Take $\frac{1}{2} < \theta < 1$ and consider first the range $n^\theta \leq j \leq n$. Using

$$\sum_{i=1}^{j} M_n(i) \leq \sum_{i=1}^{n} D_i,$$

we see that

(2.25) $$\frac{1}{n^{2\theta}} \max_{n^\theta \leq j \leq n} \sum_{i=1}^{j} M_n(i) \to 0 \quad \text{a.s.}$$

by the law of large numbers if $ED < \infty$. If we assume $ED = \infty$ and (1.5) instead, then (2.25) follows from Theorem 8.9 on page 68 of Durrett (1996), since

$$\sum_{n} P(D \geq n^{2\theta}) < \infty.$$

Using the fact that $j(j-1) \geq n^\theta(n^\theta - 1)$ for $n$ in this range, it follows under either assumption that

(2.26) $$\max_{n^\theta \leq j \leq n} \left[ 2 \sum_{i=1}^{j} M_n(i) - j(j-1) \right] \to -\infty.$$

To handle the other values of $j$, use (2.4) to write

(2.27) $$\sum_{i=1}^{j} M_n(i) \stackrel{d}{\leq} \sum_{i=1}^{j} g\left(\frac{n}{S_i}\right) \stackrel{d}{\leq} \sum_{i=1}^{j} g\left(\frac{n}{iS}\right),$$

where

$$S = \inf_{i \geq 1} \frac{S_i}{i}.$$

Note that $S > 0$ a.s. by the law of large numbers.

Assume now that $ED < \infty$. By the monotonicity of $g$,

$$\frac{1}{n} g\left(\frac{n}{iS}\right) \leq \int_{(i-1)/n}^{i/n} g\left(\frac{1}{xS}\right) dx.$$

Therefore, by (2.27),

(2.28) $$\frac{1}{n} \sum_{1 \leq i \leq n^\theta} M_n(i) \stackrel{d}{\leq} \int_{0}^{n^{\theta-1}} g\left(\frac{1}{xS}\right) dx.$$

The right-hand side of this expression tends to 0 as $n \to \infty$ provided that

$$\int_{0} g(1/x) \, dx < \infty,$$

or, equivalently, that

$$\int^\infty \frac{g(x)}{x^2}\,dx < \infty.$$

Since $g(x) = n$ for $g^{-1}(n) \leq x < g^{-1}(n+1)$, this is equivalent to

$$\sum_n \int_{g^{-1}(n)}^{g^{-1}(n+1)} \frac{n}{x^2}\,dx = \sum_n n\left[\frac{1}{g^{-1}(n)} - \frac{1}{g^{-1}(n+1)}\right] < \infty.$$

This follows by summation by parts from $ED < \infty$ and (2.3). So, the left-hand side of (2.28) tends to zero and hence

$$\max_{1 \leq j \leq n^\theta}\left[2\sum_{i=1}^{j} M_n(i) - n\right] \to -\infty.$$

Combining this with (2.26) gives (2.24) when $ED < \infty$.

Assume now that $ED = \infty$ and (1.5) holds. By (1.5) and (2.3), there is a constant $C$ so that

$$g(x) \leq \frac{Cx}{\log(x+1)}, \qquad x \geq 0.$$

Therefore, by (2.27),

$$\frac{1}{n}\sum_{1 \leq i \leq \log n} M_n(i) \stackrel{d}{\leq} C \sum_{1 \leq i \leq \log n}\left[iS\log\left(\frac{n}{iS}+1\right)\right]^{-1} \to 0.$$

It follows that

(2.29) $$\max_{1 \leq j \leq \log n}\left[2\sum_{i=1}^{j} M_n(i) - n\right] \to -\infty.$$

The same argument shows that

$$\frac{1}{n}\sum_{1 \leq i \leq n^\theta} M_n(i) \stackrel{d}{\leq} C\left[\log\left(\frac{n^{1-\theta}}{S}+1\right)\right]^{-1} \sum_{1 \leq i \leq n^\theta} \frac{1}{i},$$

which remains bounded as $n \to \infty$. However, since $ED = \infty$,

$$\frac{1}{n}\sum_{i=1}^{n} \min(\log n, D_i) \to \infty \qquad \text{a.s.}$$

by the law of large numbers. Therefore

$$\max_{\log n \leq j \leq n^\theta}\left[2\sum_{i=1}^{j} M_n(i) - \sum_{i=1}^{n}\min(j, D_i)\right] \to -\infty.$$

Combining this with (2.26) and (2.29) gives (2.24) as required. □

DEPARTMENT OF MATHEMATICS
UNIVERSITY OF SOUTHERN CALIFORNIA
LOS ANGELES, CALIFORNIA 90089-1113
USA
E-MAIL: rarratia@math.usc.edu
URL: www.math.usc.edu/~rarratia

DEPARTMENT OF MATHEMATICS
UNIVERSITY OF CALIFORNIA
LOS ANGELES, CALIFORNIA 90095-1555
USA
E-MAIL: tml@math.ucla.edu
URL: www.math.ucla.edu/~tml